\magnification 1200
\input amssym.def
\input amssym.tex
\advance\hoffset by 2truecm
\advance\hsize by 2,3 truecm
\def\makefootline{\baselineskip=52pt\line{\the\footline}}
\vsize= 23,5 true cm
\hsize= 14 true cm
\overfullrule=0mm

\font\grandsy=cmsy10 scaled \magstep2
\def\SS{{\grandsy x}}

\font\petitsy=cmsy10 scaled \magstep0
\def\SSpetit{{\petitsy x}}
\font\grandsy=cmsy10 scaled\magstep2



\def\rcases#1{\left.\, \vcenter{\normalbaselines
\ialign{$##\hfil$&\quad##\hfil\crcr#1\crcr}}\right.}

\def\build#1_#2^#3{\mathrel{\mathop{\kern 0pt#1}\limits_{#2}^{#3}}}

\headline={\ifnum\pageno > 1\hfill\tenrm\folio\hfil\else\fi}
\footline={\hfill}\pageno=1

\newcount\coefftaille \newdimen\taille
\newdimen\htstrut \newdimen\wdstrut
\newdimen\ts \newdimen\tss

\catcode`\@=11
\def\Eqalign#1{\null\,\vcenter{\openup\jot\m@th\ialign{
\strut\hfil$\displaystyle{##}$&$\displaystyle{{}##}$\hfil
&&\quad\strut\hfil$\displaystyle{##}$&$\displaystyle{{}##}$
\hfil\crcr#1\crcr}}\,}  \catcode`@12

\def\fixetaille#1{\coefftaille=#1
\htstrut=8.5pt \multiply \htstrut by \coefftaille \divide \htstrut by 1000
 \wdstrut=3.5pt \multiply \wdstrut by \coefftaille \divide \wdstrut by 1000
\taille=10pt  \multiply \taille by \coefftaille \divide \taille by 1000
\ts=\taille \multiply \ts by 70 \divide \ts by 100
 \tss=\taille \multiply \tss by 50 \divide \tss by 100
\font\tenrmp=cmr10 at \taille
\font\sevenrmp=cmr7 at \ts
\font\fivermp=cmr5 at \tss
\font\tenip=cmmi10 at \taille
\font\sevenip=cmmi7 at \ts
\font\fiveip=cmmi5 at \tss
\font\tensyp=cmsy10 at \taille
\font\sevensyp=cmsy7 at \ts
\font\fivesyp=cmsy5 at \tss
\font\tenexp=cmex10 at \taille
\font\tenitp=cmti10 at \taille
\font\tenbfp=cmbx10 at \taille
\font\tenslp=cmsl10 at \taille}

\def\fspeciale{\textfont0=\tenrmp%
\scriptfont0=\sevenrmp%
\scriptscriptfont0=\fivermp%
\textfont1=\tenip%
\scriptfont1=\sevenip%
\scriptscriptfont1=\fiveip%
\textfont2=\tensyp%
\scriptfont2=\sevensyp%
\scriptscriptfont2=\fivesyp%
\textfont3=\tenexp%
\scriptfont3=\tenexp%
\scriptscriptfont3=\tenexp%
\textfont\itfam=\tenitp%
\textfont\bffam=\tenbfp%
\textfont\slfam=\tenbfp%
\def\it{\fam\itfam\tenitp}%
\def\bf{\fam\bffam\tenbfp}%
\def\rm{\fam0\tenrmp}%
\def\sl{\fam\slfam\tenslp}%
\normalbaselineskip=12pt%
\multiply \normalbaselineskip by \coefftaille%
\divide \normalbaselineskip by 1000%
\normalbaselines%
\abovedisplayskip=10pt plus 2pt minus 7pt%
\multiply \abovedisplayskip by \coefftaille%
\divide \abovedisplayskip by 1000%
\belowdisplayskip=7pt plus 3pt minus 4pt%
\multiply \belowdisplayskip by \coefftaille%
\divide \belowdisplayskip by 1000%
\setbox\strutbox=\hbox{\vrule height\htstrut depth\wdstrut width 0pt}%
\rm}

\parindent=0pt

\def\vmid#1{\mid\!#1\!\mid}

\def\fle{\rightarrow}
\def\flep{\longrightarrow}

\null\vskip-1cm

\font\sc=cmcsc10
\font\grandsy=cmsy10 scaled\magstep2

\newdimen\emm 
\def\pmb#1{\emm=0.03em\leavevmode\setbox0=\hbox{#1}
\kern0.901\emm\raise0.434\emm\copy0\kern-\wd0
\kern-0.678\emm\raise0.975\emm\copy0\kern-\wd0
\kern-0.846\emm\raise0.782\emm\copy0\kern-\wd0
\kern-0.377\emm\raise-0.000\emm\copy0\kern-\wd0
\kern0.377\emm\raise-0.782\emm\copy0\kern-\wd0
\kern0.846\emm\raise-0.975\emm\copy0\kern-\wd0
\kern0.678\emm\raise-0.434\emm\copy0\kern-\wd0
\kern\wd0\kern-0.901\emm}

\font\tendb=msbm10
\font\sevendb=msbm7

\newfam\dbfam
\textfont\dbfam=\tendb\scriptfont\dbfam=\sevendb\scriptscriptfont\dbfam=\sevendb
\def\db{\fam\dbfam\tendb}

\def\C{{\db C }}

\def\N{{\db N }}

\def\P{{\db P }}

\def\R{{\db R }}

\def\Z{{\db Z }}

\newdimen\margeg \margeg=0pt
\def\bb#1&#2&#3&#4&#5&{\par{\parindent=0pt
    \advance\margeg by 1.1truecm\leftskip=\margeg
    {\everypar{\leftskip=\margeg}\smallbreak\noindent
    \hbox to 0pt{\hss\bf [#1]~~}{\bf #2 - }#3~; {\it #4.}\par\medskip
    #5 }
\medskip}}

\newdimen\margeg \margeg=0pt
\def\bbaa#1&#2&#3&#4&#5&{\par{\parindent=0pt
    \advance\margeg by 1.1truecm\leftskip=\margeg
    {\everypar{\leftskip=\margeg}\smallbreak\noindent
    \hbox to 0pt{\hss [#1]~~}{\pmb{\sc #2} - }#3~; {\it #4.}\par\medskip
    #5 }
\medskip}}

\newdimen\margeg \margeg=0pt
\def\bba#1&#2&#3&#4&#5&{\par{\parindent=0pt
    \advance\margeg by 1.1truecm\leftskip=\margeg
    {\everypar{\leftskip=\margeg}\smallbreak\noindent
    \hbox to 0pt{\hss [#1]~~}{{\sc #2} - }#3~; {\it #4.}\par\medskip
    #5 }
\medskip}}


\null\vskip 1cm
\centerline{\bf JACOBI ELLIPTIC CLIFFORDIAN FUNCTIONS}
\bigskip
\centerline{\bf by}
\medskip
\centerline{{\sc Guy Laville}\footnote{*}{\sevenrm  Corresponding author e-mail:
glaville@math.unicaen.fr} and {\sc Ivan Ramadanoff}\footnote{$^\dagger$}{\sevenrm
e-mail:rama@math.unicaen.fr }
}

\vskip1,5cm
{\narrower\noindent\fixetaille{800}{\fspeciale   The well-known Jacobi elliptic
functions \
$sn(z)$, $cn(z)$, $dn(z)$ \ are defined in higher dimensional spaces by the following
method. Consider the Clifford algebra of the antieuclidean vector space of dimension \ 
$2m+1$.  Let \ $x$ \ be the identity mapping on the space of scalars + vectors. The
holomorphic Cliffordian functions may be viewed roughly as generated by the powers of \
$x$, \ namely \ $x^n$,~ their derivatives, their sums, their limits (cf : $z^n$ \ for
classical holomorphic functions). In that context it is possible to define the same type of
functions as Jacobi's.}

\bigskip
\noindent\fixetaille{800}{\fspeciale\sl Keywords}~:\  {Clifford analysis,  Elliptic
functions, Jacobi functions, Holomorphic Cliffordian functions.}

\medskip
\noindent\fixetaille{800}{\fspeciale\sl AMS Classifications}~: \  {30G35, \ 33E05.}
\par}

\vskip1,5cm
\hskip-0,9cm {\bf  Introduction}
\bigskip\medskip
The theory of elliptic functions, i.e. holomorphic periodic functions of one complex
variable is well-known. The theory of periodic functions in higher dimensional spaces,
(be they real, vector or Clifford-valued) has a long history. The Dirichlet problem in a
box was the motivation of the works of P. Appell [1], [2], [3],  A. Dixon [6]. After a long
drowsiness Fueter [8] made some studies in the context of his theory of regular
quaternion-valued functions. More recently J. Ryan [15] and S. Krausshar [9]
worked with Clifford-valued functions. Here we are looking at the well-known Jacobi elliptic
functions \ $sn(z)$, $cn(z)$, $dn(z)$ \ in higher dimensional spaces, in the framework
of what we think to be the natural context~: holomorphic Cliffordian functions. The main
tool is the fundamental \ $\zeta_N$ \ functions introduced in [12].

\vskip1cm
\hskip-0,9cm {\bf \SS 1. Perequisites}
\bigskip\medskip
In this first paragraph, we will recall the notion of holomorphic
Cliffordian function, introduced in  \ [10], [11], \ and also what we could
call an elliptic Cliffordian function, \ [12]. Some basic properties of the
Cliffordian analogous of the \ $\zeta$ \ Weierstrass function will be remembered and
some ingredients as well, which we will make use further.

\bigskip\medskip
Let \ $\R_{0,2m+1}$ \ be the Clifford algebra of the real vector space \ $V$ \ of
dimension \ $2m+1$, \ provided with a quadratic form of negative signature and \
$m\in\N$.  Denote by \ $S$ \ the set of the scalars in \ $\R_{0,2m+1}$, \ which can
be identified to \ $\R$. Let \ $\{ e_i\}$, \ $i = 1,2,\ldots,2m+1$, \ be an
orthonormal basis of \ $V$ \ and set \ $e_0 = 1$. Thus, in the algebra \
$\R_{0,2m+1}$, \ the calculus rules will be generated by \ $e_ie_j + e_je_i =
-2\delta_{ij}$ \ for \ $0\leq i,j \leq 2m+1$, \ where \ $\delta_{ij}$ \ is the
Kronecker symbol.

\bigskip\medskip
A point \ $x = (x_0, x_1,\ldots, x_{2m+1})$ \ of \ $\R^{2m+2}$ \ could be
considered as an element of \ $S\oplus V$ \ and be written as \ $x = x_0 + \vec
x$, \ where \ $x_0$ \ means its scalar part and \ $\vec x =
\displaystyle\sum_{i=1}^{2m+1} e_i x_i$ \ its vector part.

\bigskip\medskip
Let \ $\Omega$ \ be an open set of \ $S\oplus V$. A function \ $f : \Omega \fle
\R_{0,2m+1}$ \ is said to be (left) holomorphic Cliffordian in \ $\Omega$ \ if and
only if~:
$$D~\Delta^m~f(x) = 0,$$

for each \ $x$ \ of \ $\Omega$. \ Here, \ $\Delta^m$ \ means the \ $m$ \ times
iterated Laplacian \ $\Delta$ \ and \ $D$ \ is the well-known operator~:
$$D = \sum_{i=0}^{2m+1} \ e_i \ {\partial\over \partial x_i}$$ 

lying on the basis of the theory of (left) monogenic functions \ ([4], see also [7]).

\bigskip\medskip
In \ [11], \ the foundations of a theory of holomorphic Cliffordian functions
were achieved, constructing a corresponding Cauchy kernel, obtaining a Cauchy
integral representation formula allowing to derive a similar to the Taylor expansion
series. Perhaps the most significant phenomenon in this theory against the theory of
monogenic functions is the fact that the function \ $x\longmapsto x^n, \ n\in\N$ \
is a holomorphic Cliffordian one. Moreover, the function \ $x\longmapsto x^{-1}$ \
is also holomorphic Cliffordian on \ $S\oplus V \setminus \{ 0\}$, so that,
restricting us only to pointwise singularities, there is no major difficulties to
obtain a similar to the Laurent expansion series for meromorphic Cliffordian
functions, [12].

\bigskip\medskip
Take now an integer \ $N$ \ in \ $\{ 1,2,\ldots, 2m+2\}$ \ and let \ $\omega_\alpha
\in S\oplus V$ \ for \ $\alpha = 1,2,\ldots,N$. Suppose always the paravectors \
$\omega_1,\ldots,\omega_N$ \ linearly independent in \ $S\oplus V$. \ For
convenience the \ $\omega_\alpha$ \ will play the role of half periods. So, a function \ $f :
S\oplus V \fle \R_{0,2m+1}$ \ is said to be \ $N$-periodic if
$$f(x+2\omega_\alpha) = f(x)$$

\medskip
for every \ $x\in S\oplus V$ \ and \ $\alpha = 1,2,\ldots,N$. Further, let us call the set \
$2\Z^N\omega = \{ 2k\omega, k\in\Z^N\}$ \ a lattice. Here we will make use of the notations~: \
$\omega = (\omega_1,\ldots,\omega_N)$, \ for a \ $N$-uple of paravectors, \ $k =
(k_1,\ldots,k_N)$ \ for a multiindex, \ $k_\alpha\in\Z$, \ and \ $k\omega =
\displaystyle\sum_{\alpha = 1}^N  \ k_\alpha \omega_\alpha$. Obviously, for a \ $N$-periodic
function, we have~:
$$f(x+2k\omega) = f(x).$$

\medskip
Recall a general theorem of elliptic Cliffordian functions, i.e. meromorphic and \
$N$-periodic~:

{\bf the theorem for the principal parts}.- \ {\it If \ $f_1$ \ and \ $f_2$ \
are two elliptic Cliffordian functions with the same pointwise poles and the same
principal parts of their Laurent expansions on the neighborhoods of their poles, then
they differ just up to an additive constant,} [12].

\bigskip\medskip
Now recall the definition of the Weierstrass \ $\zeta_N$ \ functions. First, we need to
rearrange the lattice \ $2\Z^N\omega$ \ in a countable set~: \ ${\{ w_p\}}_{p=0}^\infty$, \
where \ $w_0 = (0,0,\ldots,0)$.  Then set~:
$$\zeta_N(x) = x^{-1} + \sum_{p=1}^\infty \ \{ (x-w_p)^{-1} + \sum_{\mu = 0}^{N-1} {(w_p^{-1}
x)}^\mu \ w_p^{-1}\}.$$

\bigskip
In such a way we have a function \ $\zeta_N$~: \  $S\oplus V \setminus 2\Z^N\omega \fle
\R_{0,2m+1}$ \ for which one can show it is a holomorphic Cliffordian one and that \ $\zeta_N$
\ possesses simple poles at the vertices of the lattice. Moreover, one has also~:
$$\zeta_N(x) = x^{-1} - \sum_{k \geq [{N\over 2}]} \ \sum_{p=1}^\infty \ {(w_p^{-1}~x)}^{2k+1}
\ w_p^{-1},$$

from which it follows that \ $\zeta_N$ \ is an odd function.

\bigskip\medskip
The function \ $\zeta_N$ \ itself is not \ $N$-periodic, but it satisfies a property of
quasi-periodicity (i.e. up to a holomorphic Cliffordian polynomial). More precisely~:
$$\zeta_N (x+2\omega) - \zeta_N(x) = 2 \ \sum_{p=0}^{[{N+1\over 2}] - 1} \ {((x+\omega) \mid
\triangledown_y)^{2p}\over (2p)!} \ \zeta_N (y)\Bigl|_{y=\omega}$$

or equivalently :
$$\zeta_N(x+\omega) - \zeta_N(x-\omega) = 2 \ \sum_{p=0}^{[{N+1\over 2}] - 1} \ {(x \mid
\triangledown_y)^{2p}\over (2p)!} \ \zeta_N(y)\Bigl|_{y=\omega}.$$

\bigskip
Let us make some comments. Here we have made use of the notation~: $(x
\vmid{\triangledown_y)^{2p} \
\zeta_N(y)}_{y=\omega}$. That means we start with the usual scalar product in \ $\R^{2m+2}$ \
of the two vectors  $x = (x_0, x_1,\ldots, x_{2m+1})$ \ and the gradient \ 
$\triangledown_y =
\bigl( \displaystyle{\partial\over \partial y_0}, {\partial\over \partial y_1} ,\ldots,
{\partial\over \partial y_{2m+1}}\bigr)$ \ which is applied to the function \ $\zeta_N$ \ as a
function of the variables \ $(y_0, y_1,\ldots, y_{2m+1})$.  Then \ $(x \mid
\triangledown_y)^{2p}$ \ means an iteration \ $2p$ \ times of \ $(x\mid\triangledown_y)$ \
and the final result is obtained substituing \ $y = \omega$.

\bigskip
In the particular case when \ $N = 2m+2$ \ which will be the most natural case appearing later,
we have~:
$$\zeta_{2m+2} (x+2\omega) - \zeta_{2m+2} (x) = 2 \ \sum_{p=0}^m \ {\bigl( (x+\omega) \mid
 \triangledown_y\bigr)^{2p}\over (2p)!} \ \zeta_N(y)\Bigl|_{y=\omega}$$

Note the right-hand side is a polynomial on \ $x$  \ of degree \ $2m$ \ which will be denoted
by \ $p_{2m} (x~; \omega)$.

\bigskip\bigskip
Finally, let us write down the Laurent expansion of \ $\zeta_N$ \ in a neighborhood of the
origin. This can be done using the formula~:
$$(x \vmid{\nabla_w)^n  \ (w^{-1})}_{w = w_{p}} = (-1)^n \ n! \quad (w_p^{-1}~x)^n \
w_p^{-1},$$ where \ $w, x, w_p \in S\oplus V$ \ and \ $p, n\in\N$. So, we get :
$$\zeta_N(x) = x^{-1} - \sum_{n=N}^\infty \ \sum_{p=1}^\infty (-1)^n \ {(x \mid
\triangledown_w)^n\over n!} (w^{-1})\Bigl|_{w=w_{p}}.$$

Because of the imparity of \ $\zeta_N$, there is a more simple formula :
$$\zeta_N(x) = x^{-1} + \sum_{k \geq [{N\over 2}]} \ {1\over (2k+1)!} \ \sum_{p=1}^\infty \ (x
\mid \triangledown_w)^{2k+1} \ (w^{-1})\Bigl|_{w=w_{p}}\leqno (1)$$

Just note that in the special case \ $N = 2m+2$, \ the first sum starts from \ $k = m+1$.

\vskip1cm
\hskip-0,9cm {\bf \SS 2.  Translations operators}
\bigskip\medskip
  
Introduce the translation  operators  $E_j$ : for a fixed lattice generated by \ $N$ \
paravectors \ $2\omega_1,\ldots,2\omega_N$ \ and for an arbitrary function \ $\varphi : S
\oplus V \flep \R_{0,2m+1}$, \ set~:
$$E_j(\varphi) (x) = \varphi (x+\omega_j), \qquad j = 1,2,\ldots,N$$

The composition of the two operators \ $E_i$ \ and \ $E_j$ \ will be denoted simply as \
$E_jE_i$.  Obviously, \ $E_jE_i = E_iE_j$. Actually, the set \ $\{ E_j\}$, $j =
1,\ldots,N$
\ generates a commutative algebra of operators. This algebra is isomorphic to the polynomial
algebra of \ $N$ \ independent variables \ $\P (X_1,\ldots,X_N)$. Remark also the square of any
translation operator gives a translation on the whole period~: \ $E_j^2 (\varphi)(x) = \varphi
(x+2\omega_j)$.

\bigskip
Suppose now \  $\varphi$ \ is a \ $N$-periodic function. We could translate this fact using the
language of the translation operators, saying \ $\varphi$ \  is \ $N$-periodic if and only
if~:
$$(I-E_j^2) (\varphi) (x) = 0$$
for \ $j = 1,2,\ldots,N$.

\bigskip\medskip
But the algebraic structure of this set of operators allows us to write :
$$I - E_j^2 = (I-E_j) (I+E_j).$$

Sometime, we could look on a special translation as, for example, \ $I - E_iE_j$.

\bigskip
In this case note that :
$$(I - E_iE_j) (\varphi) (x) = \varphi (x) - \varphi (x+\omega_i+\omega_j).$$

\bigskip
The same result would be obtained if we use the identity~:
$$I - E_iE_j = (I-E_i) (I+E_j) + (I + E_i) - (I+E_j),$$
$1\leq i,j \leq N$.

\bigskip\bigskip
The next step is to understand how to write the quasi periodicity of a function as, for
example, the function \ $\zeta_N$. But it is clear that applying the operator \ $I - E_j^2$ \
on \ $\zeta_N$, \ we will get the opposite of the polynomial giving the right hand side of the
quasi periodicity. In the case \ $N = 2m+2$, \ we will have~:
$$(I - E_j^2) (\zeta_{2m+2}) (x) = -~p_{2m} (x~;~\omega_j)$$
for \ $j = 1,2,\ldots, 2m+2$.

\bigskip\bigskip
{\bf Recapitulate} : \ {\it for a \ $N$-periodic function, the operators \ $I -
E_j^2$,
$j = 1,2,\ldots,N$, \ give zero, while for a quasi-periodic function, they generate a
polynomial}. \ For \ $\zeta_N$, \ the degree of the corresponding polynomial is \ $2 \bigl(
\bigl[\displaystyle{N+1\over 2}\bigr] - 1\bigr)$.

\bigskip\medskip
The translation operators possess a very beneficial property. Remember that any holomorphic
Cliffordian polynomial could be written as a sum of monomials of the type \ $(\lambda
x)^n\lambda$, \ where \ $\lambda, x$ \ are paravectors, \  $n\in\N$. A direct observation on \
$(I - E_j) \bigl( (\lambda x)^n\lambda\bigr)$ \ shows that the last expression should be a
polynomial of degree \ $n-1$. So, applying once \ $I-E_j$ \ on a polynomial of degree $n$, \
one get a polynomial of degree \ $n-1$. Obviously, at the end of the chain, one has
$$(I - E_j) (h) = 0,$$
for  \ \  $h\in S \oplus V$, \ \ which is a polynomial of degree $0$.

\bigskip\bigskip
Let us remark that \ {\it in order to annihilate a polynomial of degree $n$, one needs to apply
\ $n+1$ \ operators of the type \ $I-E_j$, \ but not necessarly with the same $j$.} \ For
example a polynomial of second degree could be annihilated independently by \ $(I - E_1)^3$ \
or \ $(I-E_1) (I-E_2) (I-E_3)$ \ etc.

\bigskip\bigskip
When one looks at the function \ $\zeta_{2m+2}$, \ then it is true that~:
$$\prod_{j=1}^{2m+1} (I-E_j) \quad (I - E_i^2)  \ (\zeta_{2m+2}) (x) = 0,$$
\medskip
because, reading this line from the right to the left, we have \ $I - E_i^2$ \ applied to \
$\zeta_{2m+2}$ \ which generates a polynomial of degree \ $2m$, \ and then \
$\build{\prod}_j^{} (I - E_j)$ \ annihilates this polynomial. But let us write the same in the
opposite order, namely~:
$$(I - E_i^2) \ \prod_{j=1}^{2m+1} \ (I - E_j) \ (\zeta_{2m+2}) (x) = 0.$$

This can be looked as the authentic periodicity only on \ $2\omega_i$ \ of the function 

$\build{\prod}_{j=1}^{2m+1} (I-E_j) \ (\zeta_{2m+2}) (x).$

\bigskip\medskip
In such a way, we dispose with a receipt to construct periodic functions starting by a
quasi-periodic.
\bigskip
Till now we will denote for brievity the function \ $\zeta_{2m+2}$ \ by \ $\zeta$.

\vskip1,5cm
\hskip-0,9cm {\bf \SS 3.  Other ways for getting periodic functions}
\bigskip\medskip
First consider the case \ $m = 0$ \ and \ $N = 2$.  The corresponding \ $\zeta_2$ \ function
coincides with the classical Weierstrass function in \ $\C$.  Its quasi-periodicity is realized
up to a polynomial of \ $0$ \ degree~:
$$\zeta_2 (x+\omega) - \zeta_2 (x-\omega) = 2\zeta_2(\omega).$$

As usually, \ $\omega$ \ is a generic notation for the two periods \ $\omega_1$ \ and \
$\omega_2$.

\bigskip\bigskip
Could we construct a \ $(2\omega_1, 2\omega_2)$ \ periodic function having two simple poles at
\ $\alpha$ \ and \ $-\alpha$ \ saying, with opposite residues \ $k$ \ and \ $-k$~?  The answer
is yes and the construction is simple. Set~:
$$\varphi (x) = k \ \zeta_2 (x-\alpha) - k~\zeta_2 (x+\alpha).$$
Obviously \ $\varphi$ \ has the required residues and the required poles. Verify \ $\varphi$ \
is a $2\omega$-periodic function~:
$$\eqalign{&\varphi (x+\omega) = k \ \zeta_2(x-\alpha+\omega) - k \ \zeta_2(x+\alpha+\omega) = k
\ \zeta_2(x-\alpha-\omega) + k \ 2\zeta_2(\omega) \  -\cr
&- k \ \zeta_2 (x+\alpha - \omega) - k \ 2\zeta_2 (\omega) = \varphi (x-\omega).\cr}$$

Look now at the case \ $m = 1$ \ and \ $N = 4$. Here the quasi periodicity of the corresponding
\ $\zeta_4$ \ function is guaranted up to an even polynomial of second degree~:
$$\zeta_4 (x+\omega) - \zeta_4 (x-\omega) = 2\zeta_4(\omega) + (x \mid \triangledown_w)^2 \
\zeta_4 (w)\Bigl|_{w=\omega}.$$

\bigskip
If we want to construct again a 4 periodic function with two simple poles and opposite
scalar residues, we need such a method able to destroy the polynomial. One possible way is
to use  the complexified Clifford algebra \ $\R_{0,m+1} \otimes \C$, \ i.e. \  the complex
space \ $\R_{0,3} \ \oplus i \R_{0,3}$ \ and then set~:
$$\varphi (x) = k \ \zeta_4 (x-\alpha) - k \ \zeta_4 (x+\alpha) 
+ ik \ \zeta_4 (x - i\alpha) - ik \ \zeta_4 (x+i\alpha).$$

\bigskip
Thus, \ $\varphi$ \ would have the required poles at \ $\alpha$ \ and \ $-\alpha$ \ with
residues \ $k$ \ and \ $-k$, \ respectively. Those poles belong to \ $\R_{0,3}$ \ (in fact in \
$S \oplus V$).  Of course \ $\varphi$ \ inherited also two other poles. A long, but direct
computation, carried on \ $\displaystyle{1\over k} \ \varphi (x+\omega)$, \ shows that \
$\varphi$ \ is periodic. The fact that the polynomial disappears is due just to \ $i^2 = -1$.

\bigskip\bigskip
So, the complexification method gives the result,~i.e. a way to annihilated a polynomial of
second degree. There is another way, coming from iteration processes usual in numerical
analysis~: if \ $p(x)$ \ is a polynomial of degree 2 of the real variable \ $x$ \ and \
$h\in\R$, \ then~:
$$p(x+3h) - 3p (x+2h) + 3p(x+h) - p(x) = 0.$$

\bigskip\bigskip
It is not difficult to generalize the last proposition in the Cliffordian case. The result
remains true if \ $x,h\in S \oplus V$ \ and \ $p$ \ is a holomorphic Cliffordian polynomial
of degree 2.

Now, take \ $\zeta_4$ \ in \ $\R_{0,3}$. Set~:
$$\varphi (x) = \zeta_4(x) - 3\zeta_4 (x+\beta) + 3\zeta_4 (x+2\beta) - \zeta_4 (x+3\beta)$$

\medskip
with an arbitrary   $\beta\in S \oplus V$.    Denote by   $p(x~;~\omega)$   the
polynomial
$2\zeta_4 (\omega) + (x\vmid{\triangledown_w)^2 \ \zeta_4(w)}_{w=\omega}\nobreak.$ 
We are in a position to show that \ $\varphi$ \ is \ $2\omega$ \ periodic. It suffices to
form \ $\varphi (x+\omega)$ \ and to apply the quasi periodicity of \ $\zeta_4$~:
$$\eqalign{ \varphi (x+\omega)
&= \zeta_4 (x+\omega) - 3\zeta_4 (x+\beta+\omega) + 3 \zeta_4(x+2\beta+\omega) -
\zeta_4(x+3\beta+\omega ) \cr 
&= \zeta_4(x-\omega) + p(x) - 3\zeta_4(x+\beta - \omega) -
\cr 
&- 3p(x+\beta) + 3\zeta_4 (x+2\beta - \omega) + 3p(x+2\beta)\cr
&- \zeta_4 (x+3\beta - \omega) - p(x+3\beta) \cr
&= \varphi (x-\omega).\cr}$$

\bigskip
{\bf Some additional receipts}~: imagine we dispose with a 4-periodic function \
$\varphi$,
\ when \ $m=1$ \ and \ $N=4$. Set the four periods as \ $2\omega_1, 2\omega_2, 2\omega_3,
2\omega_4$.  If we want to obtain a new 4-periodic function with some of the periods
unchanged, some of them divided by half, then it suffices to add corresponding
translations. For example~:
$$\eta(x) = (I+E_1) (\varphi) (x) = \varphi(x) + \varphi (x+\omega_1)$$

would be periodic on \ $\omega_1, 2\omega_2, 2\omega_3, 2\omega_4$. The proof is obvious.
With the same initial situation~:
$$\eqalign{ \theta (x) &= (I + E_2) (I+E_3) (I+E_4)(\varphi)(x) = \varphi (x) +
\sum_{i=2}^{4} \varphi (x+\omega_i) +\cr
&+ \sum_{2\leq i < j \leq 4} \varphi (x+\omega_i + \omega_j) + \varphi (x + \sum_{i=2}^4
\omega_i)\cr}$$

would be periodic on \ $2\omega_1, \omega_2, \omega_3, \omega_4$.

\bigskip
A last remark~: remember that in the last part of \ \SSpetit 2 \ we got a method for
constructing periodic  functions from quasi-periodic in the general frame of \
$\R_{0,2m+1}$. This would be the tool we will make use systematicaly.

\vskip1,5cm
\hskip-0,9cm {\bf \SS 4.  The Jacobi elliptic Cliffordian functions}
\bigskip\medskip

The aim of this paragraph is to build in \ $\R_{0,2m+1}$ \ a system of functions which
could be viewed as the analogous of the Jacobi elliptic functions~: \ $sn, cn$ \ and \ 
$dn$.  First of all, in the complex case, which coincides with our case \ $m = 0$, \ $N =
2$, \ they are three elliptic functions whose general characteristics are~: all of them
are 2-periodic and they have the same simple poles at the same points. In the traditional
notations [16], the two main periods are \ $4K$ \ and \ $4iK'$, \ with \ $K, K'\in\R$ \
related with the evaluation of the elliptic integral under the form of Legendre.
Furthermore, they are different because \ $sn$ \ is \ $(4K, 2iK')$ \ periodic, \ $dn$ \ is
\ $(2K, 4iK')$ \  periodic and, for \ $cn$, \ the periods are submitted to a strange
perturbation~: they are \ $(4K, 2K+i2K')$, \ or equivalently \ $(2K+i2K', 4iK')$. Note
that the three vectors \ $4K, 2K+i2K'$ \ and \ $4iK'$ \ are \ $\R$-dependent.

\bigskip\bigskip
According to the end of \ \SSpetit 2, \ we are in a position to construct periodic
functions starting from \ $\zeta$~:~the abreviated notation for \ $\zeta_{2m+2}$, \
applying products of (at least)  \ $2m+1$ \ operators of the type \ $I - E_j$. \ Such a
product \ $\build{\prod}_j^{} (I - E_j)$, \ of length \ $2m+1$, \ independently if the \
$j$ \ are equal or different,  just belonging to \ $\{ 1,2,\ldots,2m+2\}$, \ will be
enough for insuring the periodicity of \ $\build{\prod}_j^{} (I-E_j) (\zeta) (x)$ \ on the
whole system of paravectors \ $2\omega_1,\ldots,2\omega_{2m+2}$.

\bigskip\bigskip
Recall that concerning the periods of \ $sn$ \ and \ $dn$, \ there is a phenomenon
consisting in a division by two along the two directions of \ $\R^2$, \ but recall also we
have a receipt allowing us to reduce by a half a given period, (see the end of \SSpetit 3).
Thus, set~:

\vskip1cm
{\sc  Definition 1.-} \ Define :
$$S_i(x) = (I+E_i) \ \prod_{j=1\atop j\not= i}^{2m+2} (I-E_j) (\zeta) (x)$$
for \ $i = 1,2,\ldots,2m+2$.

\bigskip
We claim that \ $S_i$ \ is periodic with periods \ $2\omega_1,\ldots,2\omega_{i-1},
\omega_i, 2\omega_{i+1},\ldots,2\omega_{2m+2}$. Let us verify \ $S_i$ \ is periodic on \
$2\omega_k, \ k\not= i$. For this, take
$$(I-E_k^2) (S_i)(x) = (I+E_i)\  \prod_{j\not= i}\  (I-E_j) (I-E_k^2) (\zeta) (x)$$

and let us read the last line from the right to the left~: \ $(I-E_k^2) (\zeta) (x) = -
p_{2m} (x~; \ \omega_k)$ : a polynomial which is annihilated by the product. In this case
\
$I + E_i$ \ does not play any role.

\bigskip\bigskip
However, the role of \ $I+E_i$ \ is playing when we say thet \ $S_i$ \ will be periodic on
\ $\omega_i$.  For this, we need to show that~:
$$(I-E_i) (S_i) (x) = 0.$$

And so, take :
$$(I - E_i) (S_i) (x) = \ \prod_{j\not= i} \  (I - E_j) (I - E_i^2) (\zeta) (x)$$

and this is clearly equal to zero.

\bigskip\bigskip
At this stage, we dispose with \ $2m+2$ \ analogues of \ $sn$ \ and \ $dn$. How to find an
analogue to \ $cn$~?  Arguing that the previous functions were built via products of lenght
\ $2m+2$, the most natural way is to set~:

\vskip1cm
{\sc Definition 2.-} \ Define :
$$C(x) = \prod_{j=1}^{2m+2} \ (I - E_j) (\zeta) (x).$$

\bigskip
We claim the periods of \ $C$ \ are \ ${\{ \omega_i + \omega_k\}}_{i=1}^{2m+2}$, \ where \
$k$ \ is arbitrarly fixed in  $\{ 1,\ldots,2m+2\}$.  The proof that \ $2\omega_k$ \ is
a period for \ $C$ \ is obvious. Let us show \ $\omega_i + \omega_k, \ i = 1,\ldots,2m+2$, \
$i\not= k$ \ are periods for \ $C$~:
$$\eqalign{
&(I - E_iE_k) (C) (x) = [(I-E_i) (I+E_k) + (I + E_i) - (I + E_k)] (C) (x) =\cr
&= (I - E_i) \ \prod_{j\not= k}^{} \ (I-E_j) (I - E_k^2) (\zeta) (x) + \prod_{j\not=
i}^{} \ (I - E_j) (I-E_i^2) (\zeta) (x)\cr 
&- \ \prod_{j\not= k} \ (I-E_j) (I-E_k^2) (\zeta) (x).\cr}$$

\bigskip
At each line, \ $\zeta$ \ generates a polynomial of degree \ $2m$ \ via \ $I-E_k^2$ \ or \
$I - E_i^2$ \ and then the polynomial is annihilated by a product of at least \ $2m+1$ \
operators.

\bigskip\bigskip
Remark there is no need to associate to the fixed \ $k$ \ an appropriate function, named \
$C_k$. Even if we do this, following the definition of \ $C$, \ one has~: \ $C = C_1 =
\ldots = C_{2m+2}$.

\bigskip\bigskip
In such a way, we constructed a set of elliptic Cliffordian functions \ $\{ C, S_1,\ldots,
S_{2m+2}\}$, \ whose number \ $2m+3$ \ does not suffer any change. It was clear the set of
functions has to be at least $2m+3$. \ They can not be more because we want \
$(2m+2)$~-~periodic functions, so the number of the operators in the product must be
unchanged. The only thing theoreticaly possible is to put more than one operator of the
type \ $I+E_j$ \ in the product. Look at~:
$$F(x) = (I+E_1) (I+E_2) \ \prod_{j=3}^{2m+2} \ (I-E_j) (\zeta) (x).$$

The last would not be periodic on \ $2\omega_1$, \ even on \ $\omega_1$, \ because it
remains an annihilating product of only \ $2m$ \ operators which is not sufficient for the
destruction of \ $p_{2m}~(x~;~\omega_1)$.

\bigskip\bigskip
{\it Finally, the number of \ $2m+3$ \ functions is optimal and the rules of their constructions
are rigid.}

\bigskip\medskip
{\bf Let us raise an ambiguity}. We remarked that any product \ $\build{\prod}_j^{}
(I-E_j)$,
\ of lenght \ $2m+1$, \ with the \ $j$ \ different or equal, annihilates the
quasi-periodicity polynomials on each direction of the paravectors belonging to the lattice. 
When we defined the functions \ $S_i$, \ we made use of products only of different
operators. The reason comes from the necessity to obey to the second constraint that all our
functions need to have the same poles at the same points. A study of the number and the
position of the poles will be done in the next paragraph.

\bigskip\bigskip
Come back to the case \ $m=0$, $N=2$ \ and, of course, \ $2m+3$ \ is 3. Look at the
functions \ $C, S_1, S_2$. Actually, we started a begining of description of the similarity
between \ $C, S_1, S_2$ \ and \ $cn$, \ $dn$ \ and \ $sn$, \ respectively. At this stage,
the similarity concerns only the periods. As we said, the problem of the poles will be
studied in \SSpetit 5. Anyway, it becomes and will be clear that \ $C, S_1, S_2$ \ are
nothing else then \ $ik cn (z+iK')$, \ $idn (z+iK')$ \ and \ $ksn (z+iK')$, \ respectively,
the last being written in the traditional notations, [16].

\vskip1,5cm
\hskip-0,9cm{\bf \SS 5.  General properties of the Jacobi elliptic
Cliffordian functions}
\bigskip\medskip
Come back to the definitions of \ $C, S_1,\ldots,S_{2m+2}$. If we want to explicit each of
them, we have to be patient. Each function is a sum of \ $2^{2m+2}$ \ terms which are~: \
$\zeta(x)$, \ followed by the sum of \ $E_j (\zeta) (x) = \zeta (x+\omega_j)$, \ for  \ $j
= 1,\ldots,2m+2$, \ this sum containing \ $C_{2m+2}^1$ \ terms, then we have to add the \
$E_jE_k(\zeta)(x) = \zeta(x+\omega_j+\omega_k)$, \ $1\leq j < k \leq 2m+2$, \ whose number
is \ $C_{2m+2}^2$ \ and so on, till the last term, which is \ $E_1E_2 \ldots E_{2m+2}
(\zeta)(x) = \zeta \bigl( x + \build{\sum}_{j=1}^{2m+2} \omega_j\bigr)$. In addition, we
have also to take into account the corresponding signs. For example~:
$$\leqalignno{
&C(x) = \zeta (x) - \sum_{j=1}^{2m+2} \zeta (x+\omega_j) + \sum_{1\leq i < j \leq 2m+2}
\zeta (x+\omega_i + \omega_j)\  - &(2)\cr 
&\sum_{i < j < k} \zeta
(x+\omega_i+\omega_j+\omega_k) + \cdots + (-1)^{2m+2} \zeta \bigl( x + \sum_{j=1}^{2m+2}
\omega_j\bigr). &\cr}$$

\bigskip
A first look at \ (2) \ shows an enormous set of simple poles organized in groups~: the
first group contains only $0$, \ the second all the half periods, then they are combined by
pairs, etc~$\ldots$ \ and in the final group we have only the vertex \
$\build{\sum}_{j=1}^{2m+2} \omega_j$. A geometrical description of the set of poles would
be better~: they are nothing else that all the vertices of a hyperparallelogram spanned
over \ $0, \omega_1,\ldots,\omega_{2m+2}$, \ whose number is \ $2^{2m+2}$.

\bigskip\medskip
Another observation could be done~: all the residues are \ +1 \ or \ -1 \ and, at least at
this moment, concerning \ $C$ , \ one can say that the sum of the residues appearing in the
expanded version of the definition is zero. In fact, we have an equal number of \ +1 \ and
\ -1.

\bigskip\medskip
What will be surprising is that, because of the known periods and due also to other hide
periodicities, we will find later, each Jacobi function will be uniquely determined by just
\ $2m+2$ \ poles and the knowledge of the respective residues (being among \ +1 \ and \
-1).

\bigskip\bigskip
Remark that the expanded expression of \ $C$ \ can be written as a shorter formula~:
$$C(x) = \sum_{k=0}^{2m+2}\  (-1)^k \  \ \sum \ E_{j_{1}}\ldots E_{j_{k}} (\zeta) (x),
\leqno (3)$$

where, in the second sum, we take \ $1\leq j_1~< \cdots <~j_k \leq 2m+2.$

\bigskip\bigskip
In such a way, one can say that the residue at the pole \ $\omega_{j_{1}} + \cdots +
\omega_{j_{k}}$ \ is exactly \ $(-1)^k$.

\bigskip
Once we have adopted this formalism, it is not difficult to write down expanded expressions
of \ $S_i$~:
$$S_i(x) = \sum_{k=0}^{2m+2} \ \sum \ (-1)^{\varepsilon (j_{1},\ldots,j_{k})} \
E_{j_{1}} \cdots E_{j_{k}} (\zeta) (x) \leqno(4)$$

where :
$$(-1)^{\varepsilon (j_{1},\ldots,j_{k})} = \cases{ 
(-1)^k \ &if \ $i\notin \{j_{1},\ldots, j_{k}\}$ \cr
(-1)^{k-1} \ &if \ $i\in \{ j_{1},\ldots, j_{k}\}$ \cr}$$

\bigskip
Thus, the residue at the pole \ $\omega_{j_{1}} + \cdots + \omega_{j_{k}}$ \ is \ $(-1)^k$
\ if \ $i\notin \{ j_1,\ldots,j_k\}$ \ and \ $(-1)^{k-1}$ \ if \ $i\in \{ j_1,\ldots,
j_k\}$.

\bigskip
As we can see, the signs of the residues of \ $S_i$, \ $i = 1,\ldots, 2m+2$, \ are also
well organized. For each function there is an equal number of signs  \ + \ and signs \ -.

\bigskip
An illustration of the previous calculations in the case \ $m = 0$, \ $N=2$ \ can be
summarized in the following table~:
\medskip
$$\Eqalign{
& && &&C(x) &&S_1(x) &&S_2(x)\cr \noalign{\smallskip} 
&\hbox{First group :} &&0 &&+1 &&+1 &&+1\cr
&\hbox{Second group :} &&\cases{\omega_1\cr \omega_2\cr} &&\rcases{-1\cr -1\cr}
&&\rcases{+1\cr -1\cr} &&\rcases{-1\cr +1\cr}\cr
&\hbox{Third group :} &&\omega_1+\omega_2 &&+1 &&-1 &&-1\cr\cr}$$

\bigskip
Here, they are 4 residues. The classical theory of elliptic functions says the functions
are of order 2, i.e. the knowledge of only 2 poles with opposite residues is sufficient for
determining each function at all. And this is really the case. For \ $S_1$, \ one needs to
know the residues at \ $0$ \ and \ $\omega_2$, those in \ $\omega_1$ \ and \
$\omega_1+\omega_2$ \ follow because of the period \ $\omega_1$ \ of \ $S_1$.  The same is
true for \ $S_2$~: once we know the residues at \ $0$ \ and \ $\omega_1$, \ those in \
$\omega_2$ \ and \ $\omega_1+\omega_2$ \ will be deduced by periodicity (on \ $\omega_2$). 
Concerning \ $C$, \ look at \ $C(\omega_2) = C(\omega_2 + 2\omega_1) = C(\omega_1)$ \ and
take into account the periodicity \ $C(\omega_1+\omega_2) = C(0)$, \ so only two poles \
$\{0, \omega_1\}$ \ are sufficient. Let us agree that in the sets of poles they are \
``determining" \ poles and \ ``additional" \ poles.  For \ $C, S_1, S_2$, \ the
corresponding sets of determining poles are \ $\{ 0, \omega_1\}$, \ $\{ 0, \omega_2\}$ \
and \ $\{ 0, \omega_1\}$, \ respectively.

\bigskip
Now let us prove that \ {\bf $C$ is an odd function}. By the expanded formula for \
$C, (2)$, \ we see \ $C$ \ is a sum of \ $2^{2m+2}$ \ terms of the form \ $\zeta
(x+\omega)$, provided with their signs, where here \ $\omega$ \ is a generic notation for
a half period or a sum of half periods. When we will take \ $C(-x)$, \ they will appear
terms of the type \ $\zeta (-x+\omega)$, \ which one can transform, because of the
imparity and the quasi-periodicity as~:
$$\zeta (-x+\omega) = - \zeta (x-\omega) = - \zeta (x+\omega - 2\omega) = - \zeta
(x+\omega) + p_{2m} (x~; \ \omega).$$ 
And thus we see that, putting \ $x \longmapsto -x$, \ we have the opposite of the
respective term, the sum of which will give \  $- C~(x)$, \ with a sum of polynomials.
Remark also that
$$\zeta (-x+\omega_i+\omega_j) = - \zeta (x+\omega_i+\omega_j) + p_{2m}(x~;~\omega_i) +
p_{2m}(x~;~\omega_j).$$
 
The remaining term would be :
$$\displaylines{\sum_{j=1}^{2m+2} p_{2m}(x~;~\omega_j) - \sum_{i<j} \
\bigl(p_{2m}(x~;~\omega_i) + p_{2m}(x~;~\omega_j)\bigr) \   + \cdots - \ \cr  
\bigl( p_{2m}(x~;~\omega_1) + p_{2m}(x~;~\omega_2)
+
\cdots + p_{2m}(x~;~\omega_{2m+2})\bigr).\cr}$$

\bigskip
In order to prove this expression is zero, introduce an abreviated notation for \
$p_{2m}(x~;~\omega_i)$ \ as \ $p_i$ \ for example. Our result would be derived from the
formula :
$$\sum_{j=1}^n \  p_j - \sum_{i<j} \ (p_i + p_j) + \sum_{i<j<k} \ (p_i+p_j+p_k) - \cdots +
(-1)^{n-1} \ \sum_{j=1}^n \ p_j = 0.\leqno (5)$$

Let us prove it. For convenience, put \ $p = \displaystyle\sum_{i=1}^n \ p_i$.

\vskip1cm
{\sc Lemma.-} \  For any \ $p_1,\ldots,p_n$ \ and \ $p$ \ as below :
$$\sum_{1\leq i_{1} < \cdots < i_{k} \leq n} \ p_{i_{1}} + p_{i_{2}} + \cdots +
p_{i_{k}} = C_{n-1}^{k-1} \  p.\leqno (6)$$

\bigskip\bigskip
The proof is achieved by a recurrence on \ $n\geq 2$. For \ $n=2$, \ first, \ if \ $k = 1$,
\ $p_1+p_2 = C_1^o~p$ \ and if \ $k = 2$, \ $p_1+p_2 = C_1^1~p$. Suppose \ (6) \ is
satisfied for \ $n-1$ \ and all \ $k = 1,\ldots,n-1$. Take the sum on the left hand side of
(6). It has two type of terms~:
$$p_1 + \sum_{2\leq j_{1} < \cdots < j_{k-1} \leq n} p_{j_{1}} + \cdots + p_{j_{k-1}}$$
which gives by the recurrence hypothesis \ $C_{n-2}^{k-2}~p$. The last terms are of the
type~:
$$\sum_{2\leq i_{1} < \cdots < i_{k} \leq n} \ p_{i_{1}} + \cdots + p_{i_{k}}$$
equal to \ $C_{n-2}^{k-1} \ p$. It remains to take into account that \ $C_{n-2}^{k-2} \ +
\ C_{n-2}^{k-1} = C_{n-1}^{k-1}$ \ and thus the lemma is proved.

\bigskip
The sum of the left hand side of \ (5) \ is equal to
$$C_{n-1}^o \ p~-~C_{n-1}^1 \ p \ + \cdots + \  (-1)^{n-1} \ C_{n-1}^{n-1} \ p = p \
\sum_{\ell = 0}^{n-1} \ (-1)^\ell \ C_{n-1}^\ell = 0.$$

\bigskip
This ends the proof that the Jacobi function \ $C$ \ is odd.
\bigskip
Similar procedures can be applied in order to prove that the functions \ $S_1,\ldots,
S_{2m+2}$ \ are all odd. Because of the signs, the respective formulas of combinatorics are
a little bit more complicated.

\vskip1cm
A direct consequence of the imparity of the Jacobi functions is \ {\sl the
possibility to determine} {\sl some of their zeroes}.  In fact, if \ $\varphi$ \
is a periodic function of period \ $\omega$ \ which is odd, then \
$\displaystyle{\omega\over 2}$  is a candidate for a zero of $\varphi$.  Write \
$\varphi (x+\omega) = \varphi (x)$ \ and put \ $x = - \displaystyle{\omega\over 2}$.  Thus
 $\varphi \bigl( \displaystyle{\omega\over 2}\bigr) = \varphi \bigl( -
\displaystyle{\omega\over 2}\bigr) = -
\varphi \ \bigl( \displaystyle{\omega\over 2}\bigr)$.

\bigskip
So, remember the sequences of periods of \ $C, S_1,\ldots, S_{2m+2}$ \ respectively, and,
after elimination of those half-periods which are poles, it remains that~:
$$C \bigl( {\omega_1\over 2} + {\omega_j\over 2}\bigr) = 0, \quad j = 2,3,\ldots,2m+2$$
and
$$S_j \ \bigl( {\omega_j\over 2}\bigr) = 0, \quad j = 1,2,\ldots,2m+2.$$
Moreover, playing with the known periodicities we can increase the lists of zeroes.

\vskip1cm
Let us prove~: {\sl $S_i$~becomes zero for \ $x = \displaystyle{\omega_i\over 2}, \ \omega_j
+
\displaystyle{\omega_i\over 2}, \ j = 1,2,\ldots,2m+2, \ j\not= i$, \ modulo \
$\omega_1,\ldots, \omega_{2m+2}$}. Really, we have shown already that \ $S_i \bigl(
\displaystyle{\omega_i\over 2}\bigr) = 0$. Thanks to the periodicity, we have also :  $S_i
(x+\omega_i+2\omega_j) = S_i(x)$, \ $j\not= i$, \ $j = 1,2,\ldots, 2m+2$.  Put \ $x = -
(\omega_j + \displaystyle{\omega_i\over 2})$ \ and use the imparity of \ $S_i$. \ Then \
$S_i \bigl( \omega_j + \displaystyle{\omega_i\over 2}\bigr) = - S_i \bigl( \omega_j +
\displaystyle{\omega_i\over 2}\bigr)$ \ and so \ $S_i \bigl( \omega_j +
\displaystyle{\omega_i\over 2}\bigr) = 0$. We said modulo \
$\omega_1,\ldots,\omega_{2m+2}$. Let us verify : starting from \ $S_i \bigl(
\displaystyle{\omega_i\over 2}\bigr) = 0$, \ we could add \ $\omega_j$ \ for \ $j \not= i$ \
and \ $\omega_i$.  In the two cases \ $S_i$ \ will vanish~: in the first case because we
already proved it, the second thanks to the periodicity on \ $\omega_i$. Start now from \
$S_i \bigl( \omega_j + \displaystyle{\omega_i\over 2}\bigr) = 0$.  If we add \ $\omega_j$,
$j\not= i$, \  remember \ $S_i$ \ is periodic on \ $2\omega_j$, \ so~:
$$S_i \bigl( \omega_j + {\omega_i\over 2} + \omega_j\bigr) = S_i \bigl( {\omega_i\over
2}\bigr) = 0.$$

\medskip
Finally, let us add \ $\omega_i$. \ But using now the periodicity of \ $S_i$ \ on \
$\omega_i$, \ we have again zero.

\bigskip
	In the same way, we can prove \ {\sl  $C$ \ becomes zero for \ $x = \displaystyle{1\over 2}
(\omega_1 + \omega_j)$, \ $j = 2,\ldots,2m+2$ \ and for \  $x = \displaystyle{3\over 2} \
\omega_1 + \displaystyle{\omega_j\over 2}$, \ modulo \ $\omega_1,\ldots,\omega_{2m+2}$. }

\bigskip
We do not affirm have found all the zeroes of the Jacobi functions, we just know these are
surely zeroes.

\bigskip
Now, let us reach the {\sl problem of the hide periodicities}. \ Start with the function \
$C$. Remember that, by construction, it posseses the following periods~: $2\omega_1,
\omega_1+\omega_2,\ldots,\omega_1+\omega_{2m+2}$. Moreover, \ $C$ \ is constructed by \
$2^{2m+2}$ \ poles. We claim that only \ $2m+2$ \ poles are sufficient for determining \
$C$ \ at all and they are \ $0, \omega_2,\ldots,\omega_{2m+2}$, \ those we called
determining poles. Really, look at \ $C(x+\omega_1) = C(x-\omega_1 + 2\omega_1) =
C(x-\omega_1) = C(x-\omega_1+\omega_1+\omega_j) = C(x+\omega_j)$ \ $j = 2,\ldots,2m+2$ \
and put \ $x=0$. The residue of \ $C$ \ at \ $\omega_1$ \ is determined as being the
residue at any point \ $\omega_j$, \ and we know it is \ $-1$.

\bigskip
Then, the residues at \ $\omega_1+\omega_j$ \ are the same as the residue at \ $0$ \
because of the ``official" periods \ $\omega_1+\omega_j$. Look now at the points \
$\omega_j+\omega_k$ \ with \ $j,k = 2,\ldots,2m+2$, \ $j<k$. Take \ $C(x + \omega_j +
\omega_k) = C(x+\omega_j + \omega_k + 2\omega_1) = C(x+\omega_1+\omega_j +
\omega_1+\omega_k) = C(x)$ \ and put \ $x=0$.

\bigskip
Further, we have~: $C(\omega_1+\omega_j+\omega_k) = C(\omega_k)$, \ or \ $C(\omega_j)$,
$C(\omega_i+\omega_j+\omega_k) = C(2\omega_1+\omega_i+\omega_j+\omega_k) = C(\omega_k)$, \
where \ $2\leq i < j < k \leq 2m+2$, \ and  \ $C(\omega_1+\omega_i+\omega_j+\omega_k) =
C(\omega_j+\omega_k)$. By a chain
argument we can deduce all the residues at the
vertices of the hyperparallelogram from those in \ $0, \omega_2,\ldots,\omega_{2m+2}$.

\bigskip
Concerning the hide periodicities of \ $S_i$, \ we will mention that by tedious
calculations, one can show, in the case \ $m=1$, $N=4$, \ that~:
$$\cases{
S_1 (x+\omega_2+\omega_3) = - S_1 (x+\omega_4) \cr
S_1 (x+\omega_3+\omega_4) = - S_1 (x+\omega_2) \cr
S_1 (x+\omega_4+\omega_2) = - S_1 (x+\omega_3) \cr}$$

and also
$$S_1 (x+\omega_2+\omega_3+\omega_4) = - S_1(x).$$

For \ $S_2$, \ one get :
$$\eqalign{
&S_2 (x+\omega_j+\omega_k) = S_2(x),\cr
&S_2(x+\omega_1+\omega_j+\omega_k) = - S_2(x)\cr}$$

for \ $1\leq j < k \leq 4$.
 
\bigskip\bigskip
Finally, one may say that the sets of determining poles for \ $S_i$ \ are~: \ $\{ 0,
\omega_1,\ldots, \widehat\omega_i,\ldots,\omega_{2m+2}\}$, \ where \ $\widehat{}$ \ means we
omit the term and \ $i = 1,\ldots,2m+2$.

\bigskip\bigskip
We will end this paragraph with a remark on the sum of the residues of the Jacobi
functions, under the condition to take the sum of the residues only at the determining
poles of the function. It is easily seing that, for each Jacobi function, the sum of the
residues is  $+1 + (-1) (2m+1) = - 2m$.

\bigskip\bigskip
The study of the \ $2m+3$ \ Jacobi elliptic Cliffordian functions shows that the structure
of general elliptic Cliffordian functions seems to be complicated and subtile because even
the general theorem on the sum of the residues for classical elliptic functions appears to
be a very particular case of the Cliffordian one.

\vskip1cm
\hskip-0,9cm{\bf \SS 6.  On the Laurent expansions of the Jacobi elliptic
Cliffordian functions}
\bigskip\medskip
Following formula \ (1) \ of \SSpetit 1, the Laurent expansion of \ $\zeta$ \ in a
neighborhood of the origin is~:
$$\zeta (x) = x^{-1} + \sum_{k \geq m+1} \ {1\over (2k+1)!} \ \sum_{p=1}^\infty \ (x
\mid\triangledown_w)^{2k+1} (w^{-1})\Bigl|_{w = w_{p}}\leqno (7)$$

\bigskip
As in \ [12], \ let us resort to a formal writting of \ (7) \ considering that by
definition~:
$$(x\mid\triangledown)^{2k+1} \ \bigl( \sum_{p=1}^\infty \ w_p^{-1}) \ := \
\sum_{p=1}^\infty \ (x\mid\triangledown_w)^{2k+1} (w^{-1})\Bigl|_{w = w_{p}}$$

even in the left hand side the sum \ $\displaystyle\sum_{p=1}^\infty \ w_p^{-1}$ \ is
obviously not convergent. Even more, let us introduce the notation \ $W$ \ for \
$\displaystyle\sum_{p=1}^\infty \ w_p^{-1}$.  In such a way, we have~:
$$\zeta (x) = x^{-1} + \sum_{k\geq m+1} \ {(x\mid\triangledown)^{2k+1}\over (2k+1)!} \
(W).\leqno (8)$$

\bigskip
As we already remarked in \ [12], \ for the complex case, \ i.e. \ $m = 0$, \ $N = 2$, \
(8) \ reduces to~:
$$\zeta (z) = {1\over z} + {(z\mid\triangledown)^3\over 3!} \ \bigl( \sum_{p=1}^\infty \
w_p^{-1}\bigr) \ + \cdots ,$$

which is another way to write the well-known Laurent expansion of the Weierstrass \ $\zeta$
\ function in \ $\C$~:
$$\zeta (z) = {1\over z} - z^3 \ \bigl( \sum_{p=1}^\infty \ {1\over w_p^4}\bigr) - z^5 \
\bigl( \sum_{p=1}^\infty \ {1\over w_p^6}\bigr) - \ \cdots \ .$$

\bigskip
Our aim is to get the Laurent expansions of the Jacobi elliptic Cliffordian functions \ $C,
S_1,\ldots, S_{2m+2}$ \ in a neighborhood of the origin. Look at \ (2), \SSpetit 5~: we see
that \ $C(x) = \zeta (x) + \phi (x)$, \ where in \ $\phi$ \ we have introduced the sum of \
$2^{2m+2} - 1$ \ terms containing translations of \ $\zeta$, \ i.e. \ $\zeta (x+\omega)$.

\bigskip
Consequently, \ $\phi$ \ has no pole at the origin, so \ $\phi$ \ is a holomorphic
Cliffordian function in the considered neighborhood. Moreover, \ $\phi = C-\zeta$ \ is an
odd function, so that \ $\phi (0) = 0$. Combining the Laurent expansion of \ $\zeta$ \ and
the usual Taylor expansion of \ $\phi$~:
$$\phi (x) = {(x \mid \triangledown_w)\over 1!} \ \phi (w)\Bigl|_{w = 0} \  + \
{(x\mid\triangledown_w)^3\over 3!}~\phi (w)\Bigl|_{w = 0} + \cdots $$

we can deduce :
$$\eqalign{
&C(x) = x^{-1} \ + \ {(x\mid\triangledown_w)\over 1!} \ \phi (w)\Bigl|_{w=0} \ + \
{(x\mid\triangledown_w)^3\over 3!} \ \phi(w)\Bigl|_{w=0} \ + \cdots +\cr
&\cdots + \ \sum_{k\geq m+1} \ {1\over (2k+1)!} \ \left[ (x\mid\triangledown)^{2k+1} \ (W)
+ (x\mid\triangledown_w)^{2k+1} \ \phi (w)\Bigl|_{w=0}\right].\cr}$$

\bigskip
Obviously, the same procedure can be applied to \ $S_i (x) = \zeta (x) + \psi_i(x), \
i=1,\ldots,2m+2$, \ in order to deduce the Laurent expansions of \ $S_i$.

\bigskip\bigskip
Remark also the conditions \ $\phi (0) = \psi_i(0) = 0$, \ $i = 1,\ldots,2m+2$, \ lead to a
numerous quantity of relations concerning the behavior of \ $\zeta$ \ in its half periods.
For example, \ $\phi (0) = 0$ \ means that~:
$$\sum_{k=1}^{2m+2} \ (-1)^{k-1} \ \sum_{1\leq i_{1} < \cdots < i_{k} \leq 2m+2} \ \zeta
(\omega_{i_{1}} + \cdots + \omega_{i_{k}}) = 0.$$

\bigskip\bigskip
Let us see how looks this relation in the case \ $m=0$, \ $N = 2$~:
$$\zeta_2 (\omega_1) + \zeta_2 (\omega_2) = \zeta_2 (\omega_1+\omega_2).$$

That is a right formula which admits a direct proof setting \ $x = 0$ \ in~:
$$\zeta_2 (x+\omega_1+\omega_2) = \zeta_2 (x-\omega_1-\omega_2) + 2\zeta_2 (\omega_1) +
2\zeta_2 (\omega_2)$$

and using the fact \ $\zeta_2$ \ is odd.

\vskip1cm
\centerline{\bf References}
\newdimen\margeg \margeg=0pt
\def\bba#1&#2&#3&#4&#5&{\par{\parindent=0pt
    \advance\margeg by 1.1truecm\leftskip=\margeg
    {\everypar{\leftskip=\margeg}\smallbreak\noindent
    \hbox to 0pt{\hss [#1]~~}{{\sc #2} - }#3~; {\it #4.}\par\medskip
    #5 }
\medskip}}

\vskip1cm

\bba 1&P. Appell&Sur les fonctions de trois variables rŽelles satisfaisant ˆ l'Žquation
diffŽrentielle \ $\Delta F = 0$&Acta Matematica, 4, (1884), 313-374& &

\bba 2&P. Appell&Sur quelques applications de la fonction \ $Z(x,y,z)$ \ ˆ la physique
mathŽmatique&Acta Matematica, 8, (1886), 265-294& &

\bba 3&P. Appell&Sur les fonctions harmoniques ˆ trois groupes de pŽriodes&Rendiconti del
circolo matematico di Palermo, 22, (1906), 361-370& &

\bba 4&F. Bracks, R. Delanghe, F. Sommen&Clifford analysis&Pitman, (1982)& &

\bba 5&C.A. Deavours&The quaternion calculus&Amer. Math. Monthly 80, (1973), 995-1008& &

\bba 6&A. Dixon&On the Newtonian Potential&Quaterly Journal of Mathematics 35,  
(1904), 283-296& &

\bba 7&R. Fueter&Die Funktionentheorie der Differentialgleichungen \ $\Delta u = 0$ \ und
\ $\Delta\Delta u = 0$ \ mit vier reellen Variablen&Comment. Math. Helv, (1935), 320-334&&

\bba 8&R. Fueter&†ber vierfachperiodische Functionen Monatshefte&Math. Phys. 48, (1939), 161-169
&&

\bba 9&R.S. Krausshar&Eisenstein series in Clifford analysis&Thesis, (2000)&&

\bba10&G. Laville, I. Ramadanoff&Fonctions holomorphes Cliffordiennes&C.R. Acad.
Sc. Paris, 326, s\'erie I (1998), 307-310& &

\bba11&G. Laville, I. Ramadanoff&Holomorphic Cliffordian Functions&Advances in
Applied Clifford Algebras, 8, n$^\circ$2 (1998), 323-340& &

\bba12&G. Laville, I. Ramadanoff&Elliptic Cliffordian Functions&(to appear in Complex
Variables)& &

\bba 13&H. Leutwiller&Modified quaternionic analysis in \ $\R^3$&Complex variables 20,
(1992), 19-51& &

\bba 14&H. Leutwiller&Quaternionic analysis versus its hyperbolic modification&(2000), preprint&
&

\bba 15&J. Ryan&Clifford analysis with generalized elliptic and quasi elliptic
functions&Appl. Anal., 13 (1982), 151-171& &

\bba 16&E.T. Whittaker, G.N. Watson&A Course of Modern Analysis&Cambridge
       University Press, (1996)& &

\vskip1cm
{\obeylines

Universit\'e de Caen - CNRS
FRE 2271
Laboratoire SDAD
DŽpartement de MathŽmatiques
Campus II
14032 Caen Cedex
France

\bigskip
glaville@math.unicaen.fr
rama@math.unicaen.fr

}

\end